\documentclass[11pt,reqno]{amsart}
\usepackage {amsmath,amsfonts,enumerate, epsfig, epsf,  indentfirst,  graphicx,euscript, amssymb,amsthm,amscd}

\usepackage[latin1]{inputenc}

\makeindex

 \theoremstyle{plain}

\theoremstyle{definition}

\newtheorem{theorem}{Theorem}

\newtheorem{problem}{Problem}

 \def\b{\mathbf}
\title[  Monge's Ellipsoid ]{ Historical Comments on Monge's Ellipsoid
and the Configurations of Lines of  Curvature  on Surfaces
Immersed in $\mathbb R^3$ }

\author[J. Sotomayor]{Jorge Sotomayor}


 \thanks{The  author is
fellow of  CNPq.
 This work was done under the project CNPQq/PADCT 620029/2004-8.}

\begin{document}

  \begin{abstract}
This is an essay on the historical landmarks leading to the study
of  principal
 configurations on surfaces, their structural stability
and further  generalizations. Here it is  pointed out that in the
work of Monge, 1796,  are found elements of the qualitative theory
of differential equations, founded by Poincar\'e in 1881. Two open
problems are proposed.
\end{abstract}
 \maketitle
\vskip .4cm

\section{Introduction }
\label{sec:1}

The book on differential geometry of  D. Struik   \cite {St}
contains key  references to the classical works on {\em principal
curvature lines} and their {\em umbilic } singularities due  to L.
Euler, G. Monge, C. Dupin, G. Darboux and A. Gullstrand, among
others (see also \cite{Sp} and, for additional references, also
\cite{r}). These papers
---notably that of Monge, complemented with Dupin's--- can be
connected with aspects of the {\em  qualitative theory of
differential equations} ({\em QTDE} for short) initiated by H.
Poincar\'e \cite{pome} and consolidated with the study of the {\it
structural stability and genericity} of differential equations in
the plane and on surfaces, which was made systematic from $1937$
to $1962$ due to
 the seminal works of Andronov  Pontrjagin and Peixoto (see \cite {ap, mp}).

 I established this  connection  by 1970, after being prompted by
the fortunate reading of Struik \cite{St}.   The
essay   \cite{moe}, written in the  Portuguese short story style,
 I recount how the interpretation of the historical landmarks led
 to inquire about the  principal
 configurations on smooth surfaces, their structural stability and
 bifurcations.

This paper
 discusses the  historical sources for the  work on
 the  structural stability of  principal curvature
 lines and umbilic points,  developed and further extended with the collaboration
 of C. Gutierrez   \cite{gs1, gsln, gs2} (see  also  the papers devoted to  other  differential
   equations
   of classical geometry: the  asymptotic lines \cite {a1, a2}, and  the  arithmetic, geometric
    and harmonic  mean curvature lines \cite {m, g, h}). This direction  of research
    led me
and R. Garcia  to the study of general mean curvature lines in
\cite{me}.

Here it is also   pointed out that in the work of Monge,
\cite{mo}, are found elements of the qualitative theory of
differential equations, founded by Poincar\'e in \cite{pome}.

This paper contains a reformulation of
 the essentials of \cite {moe}, updates    the  references and proposes two open problems.

\section{The Landmarks before Poincar\'e: Euler, Monge and  Dupin}
\label{sec:2}

\begin{figure}[htb]
\begin{center}
\includegraphics[height=4cm, width=6cm]{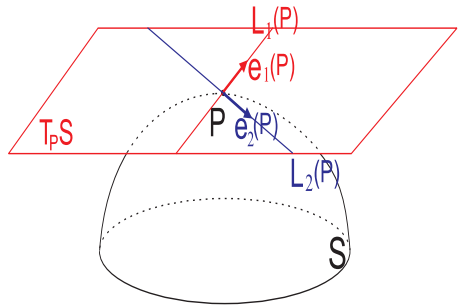}
\end{center}
\caption{Principal Directions \label{fig:1}}
\end{figure}

\noindent {\bf Leonhard Euler} (1707-1783) \cite{le}, founded of
the curvature  theory of surfaces. He defined the {\em normal
curvature} $k_n (p,L)$ on an oriented surface {\bf S} in a tangent
direction $L$ at a point $p$ as the
 curvature, at $p$,   of the planar curve of intersection of the surface with the
 plane generated by the line $L$ and the positive unit normal $N$ to the surface at $p$.
 The  {\em   principal curvatures} at $p$ are the extremal values of  $k_n (p,L)$ when
 $L$ ranges over  the tangent directions through $p$. Thus,     $k_1 (p)  = k_n (p, L_1) $
 is the {\em  minimal}   and $k_2 (p)= k_n (p, L_2) $ is the {\em maximal normal curvatures},  attained
 along the  {\em principal directions}: $L_1 (p)$,  {\em the minimal}, and  $L_2 (p)$,
 {\em  the maximal} (see Fig. \ref{fig:1}).

Euler's formula expresses the normal curvature $k_n (\theta)$
along a direction making angle $\theta$ with the minimal principal
direction  $L_1$ as  $k_n (\theta)=k_1 \cos ^2 \theta  + k_2
\sin^2 \theta$.

Euler, however,  seems to have not  considered the integral curves
of the principal line fields $L_i : p\to L_i (p), \,  i\, =1,\, 2,
\, $ and overlooked the role  of the {\em umbilic points}  at
which the principal curvatures coincide and the line fields are
undefined.

\begin{figure}[htb]
\begin{center}
\includegraphics[height=4cm, width=6cm]{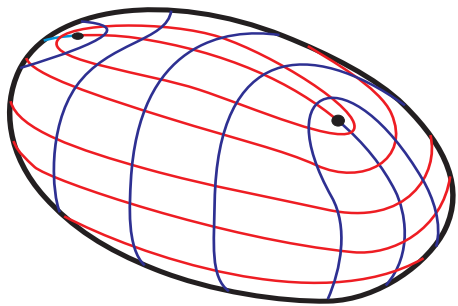}
\end{center}
\caption{Monge's Ellipsoid \label{fig:2}}
\end{figure}

\vskip 0.2cm \noindent {\bf Gaspard Monge} (1746-1818)  coined the
mathematical term {\em umbilic point} in the sense defined  above
and found the family of integral curves of the {\em principal line
fields}  $L_i ,    i\, =1,\, 2, \,  $ for the case of the triaxial
ellipsoid
$$\frac{x^2}{a^2} + \frac{y^2}{b^2} + \frac{z^2}{c^2}  -
1 = 0, \,\,\, a>b>c>0.$$ In doing this, by direct integration of
the differential equations of the principal curvature lines, circa
1779, Monge was led to the first example of a  {\em foliation with
singularities} on a surface  which ( from now on) will be called
the {\em principal configuration} of the oriented surface.  The
Ellipsoid, endowed with its principal configuration, will be
called Monge's Ellipsoid (see Fig. \ref{fig:2}).

\vskip 0.2cm

The motivation found in Monge's paper \cite {mo} is a complex
interaction of esthetic and  practical considerations and of the
explicit desire  to apply the results of his mathematical research
to real world problems. The principal configuration on the
triaxial ellipsoid appears in Monge's proposal  for the dome of
the Legislative Palace for the government of the French
Revolution, to be built over an elliptical terrain. The lines of
curvature are the guiding  curves for the workers to put the
stones. The umbilic points, from which were to hang the candle
lights, would also be the reference points below which to put the
podiums for  the speakers.

\vskip 0.2cm

 The ellipsoid depicted in Fig. \ref{fig:2} contains some of the  typical
features of the qualitative theory of differential equations
discussed briefly in {\bf a)} to {\bf d)} below:

{\bf a)  Singular Points and Separatrices.}   The umbilic points
play the role of singular points for the principal foliations,
each of them has one separatrix for each principal foliation. This
separatrix produces a connection with another umbilic point of the
same type, for which it is also a separatrix, in fact  an umbilic
separatrix connection.

\vskip 0.2cm

{\bf b)  Cycles.} The configuration  has principal cycles. In
fact, all the principal lines, except for the four umbilic
connections, are periodic. The cycles fill a cylinder or annulus,
for each foliation. This pattern is common to all classical
examples, where no surface exhibiting an isolated cycle was known.
This fact seems to be derived from the symmetry of the surfaces
considered, or from the integrability that is present in the
application of Dupin's Theorem for triply orthogonal families of
surfaces.

 As was  shown
in \cite{gs1}, these configurations an exceptional; the generic
principal cycle  for a smooth surface is  a  hyperbolic limit
cycle (see below).

\vskip 0.2cm
 {\bf c) Structural Stability (relative to quadric
surfaces).} The principal configuration remains qualitatively
unchanged under small perturbations on the coefficients of the
quadric polynomial that defines the surface.
 \vskip 0.2cm

 {\bf d) Bifurcations.}
The drastic changes in the principal configuration exhibited by
the transitions from a sphere, to an ellipsoid of revolution and
to a triaxial ellipsoid (as in  Fig. \ref{fig:2}), which after  a
very small perturbation, is a simple form of a bifurcation
phenomenon.

\vskip 0.2cm

\begin{figure}[htb]
\begin{center}
\includegraphics[height=4cm, width=6cm]{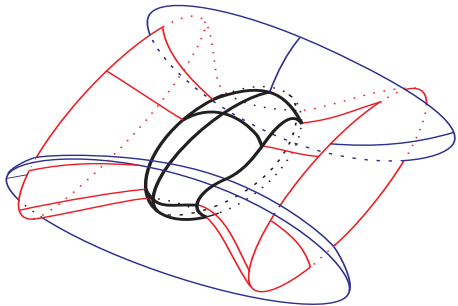}
\end{center}
\caption{Dupin's Theorem \label{fig:3}}
\end{figure}

\noindent {\bf Charles Dupin } (1784-1873) considered the surfaces
that  belong to {\em triply orthogonal surfaces,}  thus extending
considerably  those whose  principal configurations can be found
by integration. Monge's Ellipsoid belongs to the family of {\em
homofocal quadrics} (see \cite{St} and  Fig. \ref{fig:3}).

The conjunction of  Monge's analysis and
Dupin extension  provides the first global theory of integrable
principal configurations, which for quadric surfaces gives those
which are also {\em principally structurally stable } under small
perturbations of the coefficients of their quadratic  defining
equations.

\begin{theorem}\label{th:1} \cite {moe}
In the space of oriented quadrics, identified with the
nine-dimensional sphere,  those having principal structurally
stable configurations are open and dense.
\end{theorem}

\vskip 0.2cm

\noindent{\bf Historical Remark } {\em The global study of  lines
of principal curvature leading to Monge's Ellipsoid, which is
analogous of the phase portrait of a differential equation,
contains elements of Poincar\'e´s  QTDE, 85 years earlier.}

\vskip 0.2cm
 This remark seems to have been  overlooked by Monge's
scientific
 historian Ren\' e  Taton (1915-2004) in his remarkable book
\cite{ta}.

\section {Poincar\'e  and  Darboux }

The exponential role played
by {\bf Henri Poincar\'e} (1857-1912) for the {\em QTDE} as well
as for other branches of mathematics is well known and has been
discussed and analyzed in several places (see for instance
\cite{fb} and \cite{ppr}).

\begin{figure}[htb]
 \begin{center}
 \includegraphics[height=4cm, width=12cm]{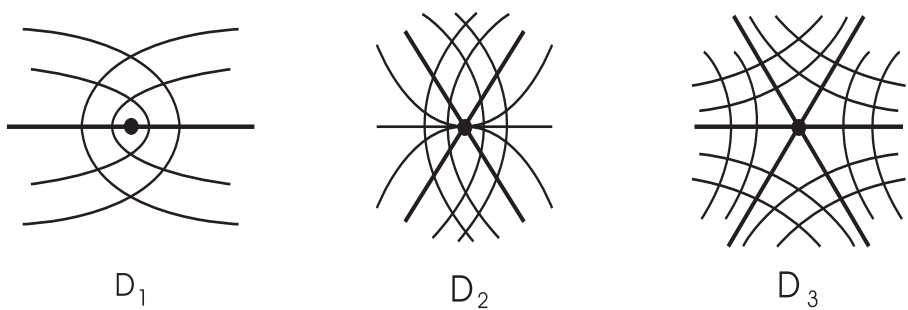}
 \caption{ Darbouxian Umbilics \label{fig:4}}
   \end{center}
 \end{figure}

\begin{figure}[htbp]
 \begin{center}
 \includegraphics[height=6cm, width=12cm]{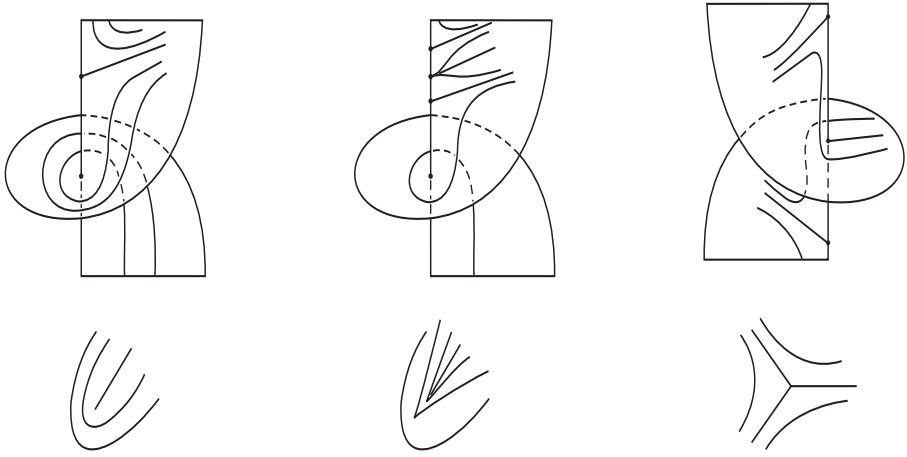}
 \caption{Lie-Cartan Resolution of Darbouxian Umbilics \label{fig:5}}
   \end{center}
 \end{figure}

Here we are concerned with his Memoires  \cite{pome}, where he
laid the foundations of the {\em QTDE}.  In this work Poincar\'e
determined the form of the solutions of planar analytic
differential equations near their  {\em foci},  {\em nodes} and
{\em saddles}. He also  studied   also properties of the solutions
around cycles and, in the case of polynomial differential
equations, also the behavior at infinity.
 \vskip 0.2cm \noindent
{\bf Gaston Darboux} (1842-1917) determined the  structure of the
lines of principal curvature near a {\em generic} umbilic point.
In his note \cite{da}, Darboux uses the theory of singularities of
Poincar\'e.  In fact, the
Darbouxian umbilics are those whose resolution 
by blowing up are saddles and nodes (see Figs. \ref{fig:4} and
\ref{fig:5}).

\vskip 0.2cm

Let $\b  p_0\in \b S$ be an umbilic point. Consider a chart
$(u,v):(\b S,\b  p_0)\to (\b R^2,\b  0)$ around it, on which the
surface has the form of the graph of a function such as
$$
\frac{k}{2}(u^2+v^2)+ \frac{a}{6}u^3+ \frac{b}{2}uv^2+
\frac{c}{6}v^3+ O[(u^2+v^2)^2].
$$

\noindent This is achieved by projecting $\b S$ onto $T\b S(\b
p_0)$ along $N(\b  p_0)$ and choosing there an orthonormal chart
$(u,v)$ on which the coefficient of the cubic term $u^2v$
vanishes.

An umbilic point is called {\it Darbouxian\/}  if, in the above
expression, the following 2 conditions {\bf T)} and {\bf D)} hold:

{\bf T)} $ b(b-a)\neq0 $

{\bf D)} either

$D_1: a/b > (c/2b)^2+2$,

$D_2: (c/2b)^2+2 >a/b >1$,  $a\neq2b$,

$ D_3: a/b <1. $

The corroboration  of the pictures in Fig.\ref{fig:4}, which
illustrate the principal configurations near Darbouxian umbilics,
has been given in \cite{gs1, gs2}; see also \cite{bf} and  Fig.
\ref{fig:5} for the Lie-Cartan resolution of a Darbouxian umbilic.

The subscripts  refer to the number of {\em umbilic separatrices,}
 which are the curves, drawn with heavy lines, tending to the
umbilic point and separating regions whose principal lines have
different patterns of approach.

\vskip 0.2cm

\section{Principal Configurations on Smooth Surfaces}\label{sec:gs}

After the seminal work of Andronov-Pontrjagin \cite{ap} on
structural stability of differential equations in the plane and
its extension to  surfaces by Peixoto \cite{mp} and in view of
the discussion on Monge's Ellipsoid formulated above, an  inquiry
into  the characterization of the oriented surfaces in $\b S$
whose principal configuration are structurally stable under small
$C^r$ perturbations, for $r\geq3$, seems unavoidable.

\vskip 0.2cm

 Call $\Sigma(a,b,c,d)$ the set of smooth compact oriented
surfaces $\b S$ which verify the following  conditions.

{\bf a)} All umbilic points are Darbouxian.
\vskip 0.1cm
 {\bf b)} All principal cycles are hyperbolic.
 This means that the corresponding   return map is
{\it hyperbolic\/}; that is,  its derivative is different from 1.
It has been shown in \cite{gs1} that hyperbolicity of a principal
cycle $\gamma$ is equivalent to the requirement  that

$$
\int_\gamma \frac{d\mathcal H} {\sqrt{\mathcal H^2-\mathcal K}}
\neq 0.
$$

\noindent where $\mathcal H = (k_1 +k_2)/2$ is the mean curvature
and $\mathcal K = k_1 k_2$ is the gaussian curvature.

\vskip 0.1cm
 {\bf c)} The limit set of every principal line is
contained in the set of umbilic points and principal cycles of $\b
S$.

The $\alpha$-(resp. $\omega$) {\it limit set\/}  of an oriented
principal line $\gamma$, defined on its maximal interval $\mathcal
I=(w_-,w_+)$ where it is parametrized by arc length $s$, is the
collection $\alpha(\gamma)$-(resp. $\omega(\gamma)$) of limit
point sequences of the form $\gamma(s_n)$, convergent in $\b S$,
with $s_n$ tending to the left (resp. right) extreme of $\mathcal
I$. The {\it limit set\/}  of $\gamma$ is the set
$\Omega=\alpha(\gamma)\cup\omega(\gamma)$.

 Examples of surfaces with {\it non
trivial recurrent\/}  principal curves, which violate condition
{\bf c} are  given in \cite{gsln, gs2}. There are no examples of
these situations in the classical geometry literature.

\vskip 0.1cm
 {\bf d)} All umbilic separatrices are separatrices of
a single umbilic point.

 Separatrices which violate {\bf d} are called {\it umbilic
connections\/}; an example can be seen in the ellipsoid of Fig.
\ref{fig:2}.

\vskip 0.2cm To make precise the  formulation  of the next
theorems, some topological notions must be defined.

 A sequence $\b S_n$ of
surfaces {\it converges in the $C^r$ sense} to a surface $\b S$
provided there is a sequence of real functions $\b f_n$ on $\b S$,
such that $\b S_n= (I+ f_n N)(\b S)$ and $f_n$ tends to $0$ in the
$C^r$ sense; that is, for every chart $(u,v)$ with inverse
parametrization $X$, $f_n\circ X$ converges to $0$, together with
the partial derivatives of order $r$, uniformly on compact parts
of the domain of  $X$.

A set $\Sigma$ of surfaces is said to be {\it open\/}  in the
$C^r$ sense if every sequence $\b S_n$ converging to $\b S$ in
$\Sigma$ in the $C^r$ sense is, for $n$ large enough, contained in
$ \Sigma$.

A set $\Sigma$ of surfaces is said to be {\it dense}  in the $C^r$
sense if, for every surface $\b S$, there is a sequence $\b S_n$
in $\Sigma$ converging to $\b S$ the $C^r$ sense.

A surface $\b S$ is said to be $C^r$-{\it principal structurally
stable\/}  if for every sequence $\b S_n$ converging to $\b S$ in
the $C^r$ sense, there is a sequence of homeomorphisms $H_n$ from
$\b S_n$ onto $\b S$, which converges to the identity of $\b S$,
such that, for $n$ big enough, $ H_n$ is  a principal equivalence
from $\b S_n$ onto $\b S$. That is, it maps $\b U_n$, the umbilic
set of $\b S_n$, onto $\b U$, the umbilic set of $\b S$, and maps
the lines of the principal foliations $\b F_{i,n}$, of $\b S_n$,
onto those of $\b F_i$, $i=1,2$, principal foliations for $\b S$.

\begin {theorem} \label{th:2} (Structural Stability of Principal
Configurations \cite{gs1, gs2})
 The set of surfaces
$\Sigma(a,b,c,d)$ is open in the $C^3$ sense and each of its
elements is $C^3$-principal structurally stable.
\end{theorem}

\begin {theorem} \label{th:3} (Density of Principal Structurally Stable
Surfaces, \cite{gsln, gs2})
 The set $\Sigma(a,b,c,d)$ is dense in the $C^2$ sense.
\end{theorem}

\vskip 0.2cm

 Extensions of these  results to
surfaces with generic singularities, algebraic surfaces, surfaces
in ${\mathbb R} ^4$ and higher dimensional manifolds  have been
achieved recently ( see, for example, \cite{gsalg, ax, mello,
gat}).

The bifurcations of umbilic points have been studied in
\cite{sell}, where references to other aspects of the bifurcations
of principal configurations are found.

\section {
Two Open Problems}\label{sec:5}

To conclude, two significant problems are stated.

\begin{problem} \label{pro:1}
Raise from 2 to 3 the differentiability  class in the density
Theorem \ref{th:3}.
\end{problem}

This remains among the most intractable
questions in this subject, involving  difficulties of Closing
Lemma type, \cite{pu}, which  also permeate other differential
equations of classical geometry, \cite{me}.

\begin{problem} \label{pro:2}
 Determine the class of principally  structurally stable  cubic
and higher degree surfaces.  Theorem \ref{th:1} deals with the
case of degree 2 ---quadric--- surfaces.
\end{problem}

Partial results, including the behavior of lines of curvature at
infinity in non-compact algebraic surfaces,  have been established
by Garcia and Sotomayor \cite{gsalg}.

\vskip 0.2cm
 \noindent{\bf Acknowledgement.\,}  I am  grateful
to C. Chicone, L. F. Mello and  R. Garcia  for their helpful
remarks, and also to and T. de Carvalho for his help in  drawing
the pictures.

\vskip .5cm
\author{\noindent Jorge Sotomayor\\Instituto de Matem\'{a}tica e Estat\'{\i}stica,\\Universidade de S\~{a}o Paulo,
\\Rua do Mat\~{a}o 1010, Cidade Universit\'{a}ria, \\CEP 05508-090, S\~{a}o Paulo, S.P., Brazil \\
\end{document}